\documentclass[fleqn]{article}
\usepackage{graphicx}
\usepackage{amsmath}
\usepackage{amsthm}
\usepackage{hyperref}
\usepackage{amssymb}
\usepackage{enumitem}
\usepackage{stmaryrd}
\usepackage{dutchcal}
\usepackage{mathtools}
\usepackage{relsize}
\usepackage{xcolor}

\newtheorem{thm}{Theorem}
\newtheorem{prop}{Proposition}
\newtheorem{lem}{Lemma}
\newtheorem{cor}{Corollary}

\newcommand\xb[1]{\ensuremath\mathbb{#1}}
\newcommand\xc[1]{\ensuremath\mathcal{#1}}
\newcommand\xf[1]{\ensuremath\mathfrak{#1}}
\newcommand\xr[1]{\ensuremath\mathrm{#1}}

\title{Absolute integral closures of commutative rings}
\author{Matth\'e van der Lee}

\begin{document}
\maketitle

\begin{abstract}
\noindent Absolute integral closures of general commutative unital rings are explored. All rings admit absolute integral closures, but in general they are not unique. Among the reduced rings with finitely many minimal prime ideals, finite products of domains are the only rings for which they are unique. Arguments using model theory suggest that the same holds for all infinite rings that are finite products of connected rings. Universal absolute integral closures, which contain every aic of a given ring, are shown to exist for certain subrings of products of domains.
\end{abstract}

\scriptsize{\textbf{Keywords:} \,Extension theory, Integral dependence, Applications of logic to commutative algebra.}

\scriptsize{{\textbf{2020 MSC:} 13B02, 13B21, 13L05.}
\normalsize

\section{Introduction}\label{sec:intro}

The absolute integral closure of an integral domain is a well-known concept. It has been much studied, and there are many applications. See \cite{H} for an overview.\\

A definition of an absolute integral closure for a general commutative ring was proposed recently by Martin Brandenburg, in a question he posted on MathOverflow regarding their uniqueness. He showed that every ring has an absolute integral closure. We will present the definition presently.\\

We work in the category $\xc{CRg}$ of commutative unital rings. Given a ring $R$, its set of minimal prime ideals is denoted by $\xr{min}(R)$, the set of regular elements by $\xr{reg}(R)$, the set of nilpotents by $\xr{nil}(R)$, and the group of units by $R^\times$. We write $\xr{Q}(R)$ for $\xr{reg}(R)^{-1}R$, the total ring of fractions of $R$. When $R$ is an integral domain, this is just the field of fractions. Recall that an $R\in\xc{CRg}$ is called \textit{connected} if $\xr{spec}(R)$ is connected topologically, that is, if $R$ is not the direct product of two nontrivial rings, that is, if $R$ has no idempotents other than the trivial ones, $0$ and $1$.\\

A ring extension $R\subseteq S$ is called \textit{tight} if non-zero ideals of $S$ contract to non-zero ideals of $R$. It suffices of course that this holds for all non-zero principal ideals of $S$. Like integrality, tightness is a transitive property: when $R\subseteq S\subseteq T$ is a chain of two tight ring extensions, $T$ is tight over $R$.\\

$S$ is \textit{absolutely integrally closed}, or \textit{ai closed} for short, if every monic $f\in S[X]$ factorizes as a product of monic linear polynomials in $S[X]$. Uniqueness of the factorization (up to a permutation of the factors) is not required.\\

If $S$ is also integral and tight over $R$, it is called an \textit{absolute integral closure} or \textit{aic} of $R$. Such $S$ are not uniquely determined up to isomorphism, as we will see. To see they exist, note that every ring $R$ has an ai closed integral extension $S$ (\cite{ST}). If we take an ideal $J$ of $S$ that is maximal wrt.\ the property that $J\cap R=0$, then $S/J$ is an aic of $R$.\\

The purpose of these notes is to investigate under what conditions on $R$ aics \textit{are} unique up to an $R$-\textit{isomorphism} (an isomorphism in the category $_R\xc{CRg}$ of $R$-algebras). For rings $R$ that have a unique aic, we will denote it by $R^{\,+}$\!, the usual notation for the absolute integral closure of $R$ when $R$ is a domain. When $R$ \textit{is} a domain, $R^{\,+}$ is its unique aic in the sense used here, see Prop.~\ref{prop:1} below.\\

Some of the basics concerning aics are set out in section \ref{sec:first}. We show that every subring of an arbitrary product of domains, over which the product is integral and tight, possesses a unique \textit{universal} aic, that is, an aic into which every aic of the ring can be embedded (Prop.~\ref{prop:2}). And we show that aics survive some forms of base change: localization and dividing out a prime ideal.\\

Section \ref{sec:redux} discusses reduced rings that have only finitely many minimal prime ideals. Such a ring is a subring of a finite product of domains, and therefore it has a universal aic $T$\!. Th.~\ref{thm:1} provides a criterion to determine when a given aic of the ring is isomorphic to $T$\!. For rings of this type, $R^\times$ exists iff $R$ is a product of domains. This is Th.~\ref{thm:2} in section \ref{sec:sample}. Actual counterexamples to aic uniqueness are also given and analyzed in that section.\\

Model theoretic methods are used in section \ref{sec:mt} to produce aics, and it is argued that infinite rings, except for finite products of domains, possess non-isomorphic aics. Prop.~\ref{prop:8} shows this, under the assumption that certain first-order theories related to the ring are well-behaved in a certain respect. The assumption is more or less self-evident, but for obvious reasons it is difficult to rigorously demonstrate that no first-order formula can exist that implies certain, infinitely many, others, to be specified in the section.

\section{First properties}\label{sec:first}

First, we establish that the new definition of absolute integral closure agrees with the existing one for domains.

\begin{prop}\label{prop:1}If $R$ is a domain, its unique aic is its $R$-plus, that is, the integral closure of $R$ in the algebraic closure $\overline{\xr{Q}(R)}$ of $\xr{Q}(R)$.
\begin{proof}
If $S$ is an aic of $R$, by integrality, lying-over holds for $S/R$, so there is a prime ideal $\xf{q}$ of $S$ with $\xf{q}\cap R=0$. As $S$ is tight over $R$, we have $\xf{q}=0$, so $S$ is a domain, and $\xr{Q}(R)\subseteq\xr{Q}(S)$ is an algebraic extension of fields. Hence $S$ can be embedded in $\overline{\xr{Q}(R)}$, and since   $S$ is ai closed and integral over $R$, its image in $\overline{\xr{Q}(R)}$ must simply coincide with $R^{\,+}$\!.
\end{proof}
\end{prop}

Next, for a finite set $\{R_a\mid a\in A\}$ of domains with product $P=\prod_{a\in A}R_a$, the product $P^{\,+}=\prod_{a\in A}R^{\,+}_a$\! is the unique aic of $P$. This is easily seen, since, as is usual with products, the product factors are completely unaware of each other, and properties only have to be verified in each component separately. It is as if one studies multiple independent rings simultaneously.\\

A case in point is the trivial ring $P=1$, which is its own aic, and is, in this sense, on par with the algebraically closed fields. The trivial ring is sometimes called the \textit{zero} ring and written as $0$, suggesting it is the initial object in the category $\xc{CRg}$, that maps to every other ring, while it is in fact the terminal object $1$ every other ring maps \textit{to}, and the name $0$ should be reserved for the true initial object, which is $\xb{Z}$ - but this aside. At any rate, being the terminal of $\xc{CRg}$, the ring $1$ is the direct product of the empty family of domains.\\

When $A$ is infinite, $P'\coloneqq P^{\,+}$ is not even integral over $P$, although every aic $S$ of $P$ still embeds in it. For let $e_a$ be the idempotent of $P$ that has $1$ as its $a$-th coordinate and zero elsewhere. Then the ideal $S_a=e_aS$ of $S$ is actually a ring. It is the product of infinitely many trivial rings and a single non-trivial one. Clearly, $S$ embeds in $S'\!=\prod_{a\in A}S_a$ via $s\mapsto(e_as)_{a\in A}$. The $S_a$ are aics of the $R_a$, so they can only be the rings $R^{\,+}_a$\!. We insert a lemma.

\begin{lem}\label{lem:1}If $R\subseteq T$ in $\xc{CRg}$ and $T$ is absolutely integrally closed, the integral closure $S$ of $R$ in $T$ is also ai closed.
\begin{proof}
If $f\in\overline{S}[X]$ is monic of degree $n$, it factors as $\prod_{i=1}^n(X-t_i)$ in $T[X]$. Then the $t_i$ are integral over the subring $R[c_{n-1},\cdots,c_{0}]$ of $S$ generated by the coefficients $c_i$ of $f$. But the latter are integral over $R$, hence so are the $t_i$. So $t_i\in S$ for all $i$, and $S$ is ai closed.
\end{proof}
\end{lem}

For infinite $A$, write $\overline{P}$ for the integral closure of $P$ in $P'$\!, which is $\{p'\in P'\mid\exists_{n\in\xb{N}}\forall_{a\in A}\exists_{f_a\in R_a[X]}(f_a\text{ monic}\,\wedge\,\xr{deg}(f_a)\le n\,\wedge\,f_a(p_a')=0)\}$. It is an aic of $P$. Indeed, $\overline{P}$ is ai closed by Lemma~\ref{lem:1}. And if $0\ne p\in\overline{P}$, take an $a\in A$ such that $e_ap\ne0$. Then $p$ is integral over $P$, so its $a$-th coordinate is integral over $R_a$, hence satisfies an equation with coefficients in $R_a$ that has a non-zero constant term $b_a$. And then $b=(b_x)_{x\in A}\in e_ap\overline{P}\cap P\subseteq p\overline{P}\cap P$ if we put $b_x=0$ for $x\in A-\{a\}$.\\

When $S$ is an aic of $P$ then, as above, $S\rightarrowtail S'\!=\prod_{a\in A}S_a\cong\prod_{a\in A}R^{\,+}_a=P'$\! is a morphism of $R$-algebras, so, since $S$ is integral over $P$, the image is contained in $\overline{P}$, and $\overline{P}$ is the universal aic of $P$. The image of $S$ may be a proper subset of $\overline{P}$. For, given $p'\in\overline{P}$, and $n\in\xb{N}$ plus $f_a\in R_a[X]$ monic of degree $\le n$ for all $a\in A$, such that $\forall_{a\in A}f_a(p_a')=0$, we can multiply each $f_a$ with a suitable power of $X$, and assume that $\xr{deg}(f_a)=n$ for all $a$. Then $f=(f_a)_{a\in A}\in P[X]$ is monic of degree $n$. But, in case every $f_a$ is separable, $f$ has $n^{|A|}$ roots in $\overline{P}$, and all we know about $S$ is that it is ai closed, so it has at least $n$ zeroes for $f$.\\

Now let $R$ be a subring of $P$ (with $A$ still infinite). The integral closure $\overline{R}$ of $R$ in $\overline{P}$ is ai closed (Lemma~\ref{lem:1}). But it is not necessarily tight over $R$. It is when $P$ is integral and tight over $R$, for then $\overline{R}=\overline{P}$ is tight over $P$, so over $R$. And, in that case, all the $e_a$ are in $\overline{R}$, for they are integral over $R$. But that does not mean that $\overline{R}=\overline{P}$ (it would for finite $A$). Since the $e_a$ are in $\overline{R}$, such $\overline{R}$ have zero-divisors. Finite products are a special case of infinite products, for we can always throw in infinitely many trivial factors.\\

We collect these observations.

\begin{prop}\label{prop:2}Let $P=\prod_{a\in A}R_a$ be a product of integral domains.
\begin{enumerate}[label=\normalfont{(}\normalfont\arabic*)]
\item If the index set $A$ is finite, $P^{\,+}$ exists, and it is the product of the $R^{\,+}_a$\!.
\item For infinite $A$, $P$ has a unique universal aic.
\item If $R$ is a subring of $P$ over which $P$ is integral and tight, then $R$ has a unique universal aic. It is the integral closure of $R$ in $\prod_{a\in A}R^{\,+}_a$\!.$\hfill\square$
\end{enumerate}
\end{prop}

Note that universal absolute integral closures, as far as they exist, are not necessarily unique up to isomorphism. Two universal aics of course embed into one another. But unless their construction is based on a reduction to a collection of domains $R$ and their $R^{\,+}$ (as in the present setup), they may have any number of roots for a given monic polynomial $f$ over the base ring. The cardinality of the zero set per $f$ must be the same in both aics, but that doesn't imply global isomorphism.\\

Returning to the general case, the following two propositions, valid for domains (\cite{H}), remain valid. The proof of the first one is immediate.

\begin{prop}\label{prop:3}If $R\subseteq S\subseteq T$ and $T$ is an aic of $R$, then $T$ is an aic of $S$.$\hfill\square$
\end{prop}

\begin{prop}\label{prop:4}If $S$ is an aic of $R$ and $F\subseteq\xr{reg}(R)$ is a monoid, $F^{-1}S$ is an aic of $F^{-1}R$. And $F\subseteq\xr{reg}(S)$. If $S=R^{\,+}$, one also has $F^{-1}S=(F^{-1}R)^+$.
\begin{proof}
The first statement is straightforward. For the second one, if $fs=0$ for an $f\in F$ and $0\ne s\in S$, by tightness $ss'\!=r\in R-\{0\}$ for some $s'$\!. But then $fr=0$, contradicting $f\in\xr{reg}(R)$. For the third statement, if $S=R^{\,+}$ and $T$ is an aic of $F^{-1}R$, then since $R\subseteq F^{-1}R$, the integral closure $Z$ of $R$ in $T$ is ai closed by Lemma~\ref{lem:1}. If $0\ne z\in Z$ and $r/f\in zT\cap F^{-1}R$ for some $f\in F$ and $0\ne r\in R$, say $r/f=zt$, then, because $t$ is integral over $F^{-1}R$, say $t^2+(x/f')t+x'/f'\!=0$ in $T$\!, with $x,x'\in R$ and $f'\in F$, to keep it simple, we find $f't\in Z$, hence $rf'\!=ztff'\in zZ\cap R$. So $Z$ is tight over $R$, and it follows that $Z\cong S$ as $R$-algebras. But since $T=F^{-1}Z$, we obtain $T\cong F^{-1}S$.
\end{proof}
\end{prop}

The next point was covered in the proof. The one after that deals with tightness and (de)localization.

\begin{cor}\label{cor:1}If the extension $R\subseteq S$ is tight, one has $\xr{reg}(R)\subseteq\xr{reg}(S)$.$\hfill\square$
\end{cor}

\begin{prop}\label{prop:5}If $R\subseteq S$, $F\subseteq\xr{reg}(R)$ is a monoid, and $S$ is tight over $R$, then $F^{-1}S$ is tight over $F^{-1}R$. Conversely, if $F^{-1}R\subseteq F^{-1}S$ is tight and the elements of $F$ are also regular in $S$, then $R\subseteq S$ is tight.
\begin{proof}
When $0\ne\frac{z}{f}\in F^{-1}S$ with $z\in S$ and $f\in F$, then $zs=r\in R-\{0\}$ for some $s\in S$. If $\frac{zs}{f}=0$ in $F^{-1}S$, then $f'r=0$ in $S$ for some $f'\in F$, contradicting regularity of $f'$\!. For the converse, for every $0\ne z\in S$ there exist $r\in R$, $s\in S$ and $f,f'\in F$ such that $0\ne\frac{z}{1}\frac{s}{f}=\frac{r}{f'}$ in $F^{-1}S$. Hence $0=f''(zsf'-rf)$ in $S$ for an $f''\in F$, so, because $f''\in\xr{reg}(S)$, $zsf'\!=rf\in R$. And $zsf'\ne0$, for else $\frac{zs}{f}$ would be zero in $F^{-1}S$.
\end{proof}
\end{prop}

If $J$ is an ideal of an aic $S$ of $R$, then $S/J$ can of course not be expected to be tight over $R/J\cap R$. Still, this property \textit{does} hold in case $J$ is a prime ideal.

\begin{lem}\label{lem:2}If $R\to S$ is integral, $\xf{q}\in\xr{spec}(S)$ and $J\supseteq\xf{q}$ is an ideal of $S$ with $J\cap R=\xf{q}\cap R$, then $J=\xf{q}$.
\begin{proof}
Let $\xf{p}=\xf{q}\cap R$. Then $\xf{q}\cap(R-\xf{p}) = \varnothing$, so $\xf{q}S_\xf{p}\in\xr{spec}(S_\xf{p})$. It is readily seen that $(\xf{q}S_\xf{p})\cap R_\xf{p}=\xf{p}R_\xf{p}$. But $R_\xf{p}\to S_\xf{p}$ is integral, and $\xf{p}R_\xf{p}\in\xr{max}(R_\xf{p})$, so $\xf{q}S_\xf{p}$ is a maximal ideal of $S_\xf{p}$ (for $R_\xf{p}/\xf{p}R_\xf{p}\subseteq S_\xf{p}/\xf{q}S_\xf{p}$ is also integral). But $\xf{q}S_\xf{p}\subseteq JS_\xf{p}\subsetneq S_\xf{p}$ (for $J\cap(R-\xf{p})$ is also $\varnothing$). Hence $JS_\xf{p}=\xf{q}S_\xf{p}$, and this implies that $J=\xf{q}$.
\end{proof}
\end{lem}

Now if $\xf{p}$ is a minimal prime ideal of $R$, there is a $\xf{q}\in\xr{spec}(S)$ lying over it. Since, by Lemma~\ref{lem:2}, there are no inclusion relations between the $S$-primes lying over $\xf{p}$, this $\xf{q}$ must be a minimal prime of $S$. Then, by Lemma~\ref{lem:2}, $S/\xf{q}$ is tight over $R/\xf{p}$, hence an aic of $R/\xf{p}$, and hence $S/\xf{q}=(R/\xf{p})^{+}$\!. This goes at least \textit{some} way towards aic uniqueness, as the $R/\xf{p}$ are the domains closest to $R$. So:

\begin{cor}\label{cor:2}If $S$ is an aic of $R$ and $\xf{q}\in\xr{spec}(S)$, then $S/\xf{q}=(R/\xf{q}\cap R)^{+}$\!.$\hfill\square$
\end{cor}

If in Lemma~\ref{lem:2} one only assumes that $\xf{q}$ is \textit{primary}, or even merely that $\sqrt{\xf{q}}\in\xr{spec}(S)$, the conclusion reads $\sqrt{J}=\sqrt{\xf{q}}$. Finally, we note that M. Artin's beautiful proof (\cite{H}, Thm.\ 2.6) of a result we quote here goes through verbatim.

\begin{prop}\label{prop:6}If $S$ is an aic of $R$ and $\xf{p}$ and $\xf{q}$ are prime ideals of $S$, then $\xf{p}+\xf{q}\in\xr{spec}(S)\cup\{S\}$.$\hfill\square$
\end{prop}

\section{Reduced rings with finite minimal spectrum}\label{sec:redux}

As an illustration of Prop.~\ref{prop:2}, we consider reduced rings $R$ for which $\xr{min}(R)$ is a finite set, and we discuss their (unique) universal aic.\\

The trivial ring $R=1$ is the case $\xr{min}(R)=\varnothing$, and integral domains are just the special case where $\xr{min}(R)$ is a one-point space. For if $\xr{min}(R)=\{\xf{p}\}$, then $\xf{p}=\bigcap\xr{min}(R)=\bigcap\xr{spec}(R)=\xr{nil}(R)=0$. This class of rings also includes all noetherian reduced rings.\\

Let $R\subseteq S$ be a given tight integral extension, with $S$ an ai closed ring. Then we have $\xr{nil}(S)\cap R=\xr{nil}(R)=0$, so by tightness $S$ must also be reduced.\\

For $\xf{p}\in\xr{min}(R)$, take a minimal prime ideal $\tilde{\xf{p}}$ of $S$ lying over it. Then, if $I=\bigcap_{\xf{p}\in\xr{min}(R)}\tilde{\xf{p}}$, we have $I\cap R=\bigcap\xr{min}(R)=\xr{nil}(R)=0$, hence $I=0$. So if $\xf{Q}\in\xr{min}(S)$, since $\xr{min}(R)$ is finite, the product $\prod_{\xf{p}\in\xr{min}(R)}\tilde{\xf{p}}$ exists and is contained in $I=0\subseteq\xf{Q}$. Thus  $\xf{Q}$ must be one of the $\tilde{\xf{p}}$. Therefore, $S$ is also "semiglobal", that is to say, $\xr{min}(S)=\{\tilde{\xf{p}}\mid\xf{p}\in\xr{min}(R)\}$ is finite.\\

Let $K$ and $L$ be the total rings of fractions of $R$ and $S$, respectively. So $K=\xr{reg}(R)^{-1}R$. The prime ideals of $K$ are of the form $\xf{p}K$, where $\xf{p}\in\xr{spec}(R)$ with $\xf{p}\cap\xr{reg}(R)=\varnothing$, i.e. $\xf{p}$ consists of zero divisors of $R$. These $\xf{p}$ are precisely the minimal prime ideals of $R$. For if $\xf{p}\in\xr{min}(R)$ and $r\in\xf{p}$, then $r\in\xf{p}R_\xf{p}$. Being a localisation of a reduced ring, $R_\xf{p}$ is again reduced. So $0=\xr{nil}(R_\xf{p})=\bigcap\xr{spec}(R_\xf{p})=\xf{p}R_\xf{p}$, for $\xf{p}R_\xf{p}$ is the only prime ideal of $R_\xf{p}$. (Note that $R_\xf{p}$ therefore must be a field.) So $r=0$ in $R_\xf{p}$, and hence there is an $r'\in R-\xf{p}$ for which $rr'\!=0$ in $R$. So $r$ is a zero divisor of $R$. Conversely, if $rr'\!=0$ and $r'\ne0$, then there is a $\xf{p}\in\xr{min}(R)$ with $r'\notin\xf{p}$, because $\bigcap\xr{min}(R)=0$. Therefore, $r\in\xf{p}$. So if a prime $\xf{q}$ of $R$ contains only zero divisors, it is contained in $\bigcup\xr{min}(R)$. But this a finite union, and so, by the prime avoidance lemma, the ideal $\xf{q}$ is a minimal prime of $R$.\\

Thus $\xr{spec}(K)=\{ \xf{p}K\mid\xf{p}\in\xr{min}(R)\}$. And we have $R\subseteq K$. Note that $\xf{p}K\cap R=\xf{p}$ when $\xf{p}$ is a minimal prime. Indeed, if $r=u/v$ in $K$ with $u\in\xf{p}$ and $v\in R$ regular, then there is a $w\in\xr{reg}(R)$ with $w(vr-u)=0$ in $R$, so $vr=u\in\xf{p}$. If $v\in\xf{p}$, then by the above $v$ is a zero divisor of $R$, contradiction. So $r\in\xf{p}$. As a result, $R/\xf{p}\subseteq K/\xf{p}K$.\\

It follows that $K$ is a zero-dimensional reduced ring, that is, a von Neumann regular ring. For if $\xf{p}$ and $\xf{q}$ are minimals of $R$ with $\xf{p}K\subseteq\xf{q}K$, then $\xf{p}\subseteq\xf{q}K\cap R=\xf{q}$, hence $\xf{p}=\xf{q}$. So if $\xf{p}\in\xr{min}(R)$, we have $\xf{p}K\in\xr{max}(K)$, and $K/\xf{p}K$ is a field. It is in fact the quotient field of $R/\xf{p}$. The ring $L$ is also VNR, and it contains $S$ as a subring.\\

$R\subseteq S\subseteq L=\xr{reg}(S)^{-1}S$, and $\xr{reg}(R)\subseteq\xr{reg}(S)$ in view of Cor.~\ref{cor:1}. Thus $K=\xr{reg}(R)^{-1}R\subseteq\xr{reg}(R)^{-1}S$ (since localizations are flat) $\subseteq\xr{reg}(S)^{-1}S=L$. Hence $K/\xf{p}K$ is a subfield of $L/\tilde{\xf{p}}L$. And the extension $K/\xf{p}K\subseteq L/\tilde{\xf{p}}L$ is algebraic, because $S$ is integral over $R$.\\

The natural map $K\to\prod_{\xf{p}\in\xr{min}(R)}K/\xf{p}K$ is injective, for the kernel is the intersection of all prime ideals of $K$, and $K$ is reduced. As $\xr{dim}(K)=0$, by the CRT this is actually an isomorphism, and $K$ is a finite product of fields. It is easy to see that in fact $\xr{Q}(R/\xf{p})\cong K/\xf{p}K\cong K_{\xf{p}K}\cong R_\xf{p}$.\\

The extension $R\rightarrowtail P\coloneqq\prod_{\,\xf{p}\in\xr{min}(R)}R/\xf{p}$ is integral and tight. Indeed, by the finiteness of the minimal spectrum, $\prod_{\,\xf{q}\in\xr{min}(R)-\{\xf{p}\}}\xf{q}$ exists. It is not contained in $\xf{p}$, and, therefore, $\bigcap_{\,\xf{q}\in\xr{min}(R)-\{\xf{p}\}}\xf{q}\nsubseteq\xf{p}$, for every $\xf{p}$. So if $0\ne p=(\overline{r}_\xf{p})_\xf{p}\in P$, with $\overline{r}_\xf{p}\ne0$, say, with lift $r_\xf{p}\in R-\xf{p}$, and $c$ is in $\bigcap_{\,\xf{q}\in\xr{min}(R)-\{\xf{p}\}}\xf{q}-\xf{p}$, then $cr_\xf{p}\in pP\cap(R-\{0\})$.\\

Take $A=\xr{min}(R)$ and put $R_a=R/a$ for $a\in A$. Then, by (1) of Prop.~\ref{prop:2}, $T=P^+=\prod_{\xf{p}\in\xr{min}(R)}(R/\xf{p})^+$ is the unique aic of $P$, hence, by (3) of Prop.~\ref{prop:2}, the universal aic of $R$. And $(R/\xf{p})^+=S/\tilde{\xf{p}}$ by Cor.~\ref{cor:2}. $T$ is also a subring of $\overline{K}\coloneqq\prod_{\xf{p}\in\xr{min}(R)}C_\xf{p}$, where $C_\xf{p}$ is the algebraic closure of the field $K/\xf{p}K$ (and of $L/\tilde{\xf{p}}L$). This $\overline{K}$ may be regarded as the "algebraic closure" of $R$ (or, equally, of $K$, $S$ or $L$). Clearly, $T$ is also the integral closure of $R$ in $\overline{K}$.\\

The tightness of $T$ over $R$ can also be seen directly. Let $0\ne t\in T$\!, and, for $\xf{p}\in\xr{min}(R)$, denote the $\xf{p}$-th coordinate of $t$ by $t_\xf{p}$, and pick a $\xf{p}$ with $t_\xf{p}\in C_\xf{p}$ nonzero. Take a monic $f\in R[X]$ with $f(t)=0$ in $T$\!. If the constant term $f(0)$ is in $\xf{p}$, it becomes zero in $K/\xf{p}K$, hence in $C_\xf{p}$. As $t_\xf{p}$ is a root of the image of $f$ in $C_\xf{p}[X]$ and $C_\xf{p}$ is a field, $t_\xf{p}$ is also a root of (the image of) $g=(f-f(0))/X\in R[X]$. We then replace $f$ by $g$. If the new $f(0)$ is in $\xf{p}$ again, repeat the process until $f(0)\notin\xf{p}$. With $c\in R$ as above, which is in all minimals of $R$ except $\xf{p}$, put $h\coloneqq cf$. Then $h(t_\xf{p})=0$, and for $\xf{p}\ne\xf{q}\in\xr{min}(R)$, we have $c=0$ in $K/\xf{q}K$, so $h=0$ in $C_\xf{q}[X]$. Hence $h(t)=0$ in $T$\!. But $h(0)=cf(0)\notin\xf{p}$, and therefore $h(0)$ is a nonzero element of $tT$ that is in $R$.\\

The image of $S$ in $\overline{K}$ under the composition map $S\subseteq L\rightarrowtail\prod_{\xf{p}\in\xr{min}(R)}L/\tilde{\xf{p}}L\subseteq\overline{L}=\overline{K}$ is integral over $R$, so it is contained in $T$\!, and this confirms the universality of the absolute integral closure $T$\!.\\

$T$ has $2^{|\xr{min}(R)|}$ idempotents. Denote the one with $e_\xf{p}=1$ and $0$ elsewhere by $e_{(\xf{p})}$. Then every idempotent is the sum of the elements of a subset of $\{e_{(\xf{p})}\mid\xf{p}\in\xr{min}(R)\}$. And $e_{(\xf{p})}e_{(\xf{q})}=0$ when $\xf{p}\ne\xf{q}$. Therefore, the $e_{(\xf{p})}$ form a fundamental system of orthogonal idempotents.\\

Then $S=T$ iff the $e_{(\xf{p})}$ are all in $S$. For if they are, and $t\in T$\!, fix a $\xf{p}$, and let $f\in R[X]$ be monic with $f(t)=0$. Then $f(t_\xf{p})=0$ in $C_\xf{p}$. Write $f=\prod_{1\leq i\leq n}(X-s_i)$ in $S[X]$, so that we have $f=\prod_{1\leq i\leq n}(X-(s_i)_{\xf{p}})$ in $C_\xf{p}[X]$. But then $t_\xf{p}=(s_i)_{\xf{p}}$ for some $i$, because $C_\xf{p}$ is a field. Put $s=s_i$. Since $e_{(\xf{p})}$ is in $S$, so is $e_{(\xf{p})}s$. Its $\xf{p}$-th component is $t_\xf{p}$, and it has zero for the other $\xf{q}$. So $e_{(\xf{p})}s$ is actually equal to $e_{(\xf{p})}t$. The sum of the $e_{(\xf{p})}$, taken over the $\xf{p}\in\xr{min}(R)$, is equal to $1$. So $t$, which is therefore the sum of the $e_{(\xf{p})}t$, must be in $S$.\\

To sum up:

\begin{thm}\label{thm:1}With the hypotheses and notation of the current section, one has:
\begin{enumerate}[label=\normalfont{(}\normalfont\arabic*)]
\item $T=\prod_{\xf{p}\in\xr{min}(R)}(R/\xf{p})^+$ is the universal aic of $R$.
\item An aic of $R$ is isomorphic to $T$ if and only if it contains precisely $2^{|\xr{min}(R)|}$ idempotents, that is, iff it has $|\xr{min}(R)|$ connected components.$\hfill\square$
\end{enumerate}
\end{thm}

Recall that when $R$ is a domain, so is $S$. Hence $S$ contains $2^{|\xr{min}(R)|}$ idemps, namely just the trivial ones, so this confirms once again that domains have a unique aic. (And the proposition is also valid for $R=1$, in which case $T=R$.)\\

With regard to generalizing the above construction, we note that for reduced $R$, the total ring of fractions $K=Q(R)$ is VNR iff every ideal of $R$ contained in $\bigcup\xr{min}(R)$ (that is, every ideal consisting entirely of zero divisors) is contained in some minimal prime of $R$. For such $R$, the minimal prime spectrum $\xr{min}(R)$ is automatically compact wrt.\ the topology induced by the Zariski topology on $\xr{spec}(R)$ - although that in inself is not enough for $K$ to be VNR (cf. \cite{M}, Prop.~1.15, or various places in \cite{T}, such as Lemma 3.2 and Th.~3.4).\\

But a fundamental difficulty with the construction when $\xr{min}(R)$ is infinite is that $T$ no longer needs to be tight over $R$. For, look at the ideal $I$ generated by the idempotent $e_{(\xf{p})}\in T$ for a $\xf{p}\in\xr{min}(R)$. It consists of the $t\in T$ with $t_{\xf{p}}\in T_{\xf{p}}$, that is, $t_{\xf{p}}\in C_{\xf{p}}$ is integral over $R$, and $t_{\xf{q}}=0$ in $C_{\xf{q}}$ for the remaining $\xf{q}\in\xr{min}(R)$. For $t$ to be in $R$ and non-zero, you would need $t_{\xf{p}}$ to be in $R$ and in every $\xf{q}\in\xr{min}(R)$ except in $\xf{p}$. And this must hold for each minimal $\xf{p}$ of $R$. This would appear to preclude any substantial generalization of the case that $\xr{min}(R)$ is finite.

\section{Sample rings having non-unique aics}\label{sec:sample}

We continue the notation of the previous section. If $k$ is a field, then the ring $R=k[X,Y]/(XY)$ has non-isomorphic aics.\\

In this case, we have $\xr{min}(R)=\{\xf{p},\xf{q}\}$, with $\xf{p}=(X)$ and $\xf{q}=(Y)$, if we denote the images of $X$ and $Y$ in $R$ by the same symbols. And $K/\xf{p}K=k(Y)$, with algebraic closure $C_\xf{p}$. The image of $R$ in $C_\xf{p}$ is the polynomial ring $k[Y]$. Let $T_\xf{p}$ be the integral closure of $k[Y]$ in $C_\xf{p}$, and $T_\xf{q}=k[X]^+\subseteq C_\xf{q}\cong_kC_\xf{p}$ (under $X\mapsto Y$). We then have $T=T_\xf{p}\times T_\xf{q}$ in $C_\xf{p}\times C_\xf{q}=\overline{K}$.\\

If $\overline{k}$ is the algebraic closure of $k$, there are ring homomorphisms $T_\xf{p}\xrightarrow{\varphi}\overline{k}
\xleftarrow{\psi}T_\xf{q}$ extending the maps $\varphi:k[Y]\to \overline{k}\gets k[X]:\psi$ defined by $Y\mapsto0\mapsfrom X$. Indeed, take $\xf{M}\in\xr{max}(T_\xf{p})$ lying over $(Y)\in\xr{max}(k[Y])$. Then $k[Y]/(Y)\subseteq T_\xf{p}/\xf{M}$ is an algebraic extension of fields, so $\varphi:k[Y]\to k[Y]/(Y)=k\subseteq\overline{k}$ can be extended to $\varphi:T_\xf{p}\to T_\xf{p}/\xf{M}\to\overline{k}$. Likewise for $T_\xf{q}$ and $\psi$.\\

Now let $T_0=\{(\alpha,\beta)\in T=T_\xf{p}\times T_\xf{q}\mid\varphi(\alpha)=\psi(\beta)\}$ be the pullback. It contains the image of $R$ in $T$\!. For every $r\in R$ is of the form $\lambda+Xf(X)+Yg(Y)$ with $\lambda\in k$ and $f(Z),g(Z)\in k[Z]$. This maps to $(\lambda+Yg(Y),\lambda+Xf(X))$ in $R/\xf{p}\times R/\xf{q}=k[Y]\times k[X]\subseteq T$\!, and thence to $(\lambda,\lambda)$ in $\overline{k}\times\overline{k}$ under $\varphi\times\psi$, since $\varphi(Y)=0=\psi(X)$. And $T_0$ is integral and tight over $R$. For let $t_0=(\alpha,\beta)\in T_0-R$. Say $\alpha\ne0$. Then $Yt_0=(Y\alpha,0)$ is in $T_0-\{0\}$, and we have $f(\alpha)=0$ for some monic $f=f(Z)\in k[Y][Z]$ with $f(0)\in k[Y]$ non-zero, for example for the minimal polynomial $f$ of $\alpha$ over $k(Y)$. As $k[Y][Z]\subseteq R[Z]$, we may view $f$ as a polynomial over $R$. Then $g\coloneqq Y^{\xr{deg}(f)}f(Y^{-1}Z)\in R[Z]$ has $g(Y\alpha)=0$ in $T_\xf{p}$ and $g(0)$ = $Y^{\xr{deg}(f)}f(0)\ne0$. But $Y=0$ in $T_\xf{q}$, so $g(0)$ vanishes in $T_\xf{q}$, hence $g(Yt_0)=0$ in $T_0$. This yields $0\neg(0)\in Yt_0T_0\cap R\subseteq t_0T_0\cap R$.\\

Finally, $T_0$ is ai closed. For let $f\in T_0[Z]$ be monic, of degree $n$, say. Then $f=\prod_{1\leq i\leq n}(Z-(\alpha_i,\beta_i))$ in $T[Z]$ for suitable $\alpha_i\in T_\xf{p}$ and $\beta_i\in T_\xf{q}$, as $T$ is absolutely integrally closed. But then $f=\prod_{1\leq i\leq n}(Z-(\alpha_i,\beta_{\pi(i)}))$ in $T[Z]=(T_\xf{p}\times T_\xf{q})[Z]$ for every permutation $\pi$ of the indices $1,\cdots,n$ (!). For $i<n$, let $(\eta_i,\vartheta_i)\in T_0$ be the coefficient of $Z^i$ in $f$. So $(\eta_0,\vartheta_0)=(-1)^n\prod_{1\leq i\leq n}(\alpha_i,\beta_i)$, and so on, up to $(\eta_{n-1},\vartheta_{n-1})=-\sum_{1\leq i\leq n}(\alpha_i,\beta_i)$. Then we have $(-1)^n\prod_{1\leq i\leq n}\varphi(\alpha_i)=\varphi(\eta_0)=\psi(\vartheta_0)=(-1)^n\prod_{1\leq i\leq n}\psi(\beta_i)$, and so forth, up to $-\sum_{1\leq i\leq n}\varphi(\alpha_i)=\varphi(\eta_{n-1})=\psi(\vartheta_{n-1})=-\sum_{1\leq i\leq n}\psi(\beta_i)$. It follows that in the ring $\overline{k}[Z]$ the equality $\prod_{1\leq i\leq n}(Z-\varphi(\alpha_i))=\prod_{1\leq i\leq n}(Z-\psi(\beta_i))$ holds, and hence $\varphi(\alpha_i)=\psi(\beta_{\pi(i)})$ for all $i$, for some permutation $\pi$.\\

So $T_0$ is an aic of $R$. But since $e_{(\xf{p})}=(1,0)$ is not in $T_0$, the rings $T_0$ and $T$ cannot be isomorphic. One is connected, while the other is not. Since $2^{|\xr{min}(R)|}=4$ here, it is clear there cannot be any other aics.\\

This generalizes easily.

\begin{thm}\label{thm:2}For $R$ reduced with $\xr{min}(R)$ finite, the following are equivalent.
\begin{enumerate}[label=\normalfont{(}\normalfont\arabic*)]
\item $R^{\,+}$\! exists, i.e., the absolute integral closure of $R$ is uniquely determined.
\item $R$ is a finite direct product of integral domains.
\end{enumerate}
\begin{proof}
One direction is given by (1) of Prop.~\ref{prop:2}. For the other, $R$ cannot be an infinite product $\prod_{a\in A}R_a$ of non-trivial rings, for if $\xf{p}_a\in\xr{min}(R_a)$, we would have $\xf{p}_a\times\prod_{a'\in A-\{a\}}R_{a'}\in\xr{min}(R)$. Since for finite products we can simply consider the individual factors individually, we may as well assume $R$ is connected. If $|\xr{min}(R)|\le1$, i.e.\ if $R$ is either a domain or the empty product $1$ of domains, we are done. In all other events, $R=R/\bigcap\xr{min}(R)
\rightarrowtail\prod_{\xf{p}\in\xr{min}(R)}R/\xf{p}$, the natural map, can't be an isomorphism, so by the (contraposition of the) CRT there are minimals $\xf{p}\ne\xf{q}$ that are not comaximal. Then $\xf{p}\cup\xf{q}\subseteq\xf{m}$ for some $\xf{m}$ in $\xr{max}(R)$, and the maps $\varphi:R/\xf{p}\to R/\xf{m}\eqqcolon k\to\overline{k}\gets R/\xf{m}\gets R/\xf{q}:\psi$ can be extended to $T_\xf{p}\xrightarrow{\varphi}\overline{k}\xleftarrow{\psi}T_\xf{q}$, where $T=T_\xf{p}\times T_\xf{q}\times\prod_{\xf{r}\in\xr{min}(R)-\{\xf{p},\xf{q}\}}T_\xf{r}$ in $\prod_{\xf{r}\in\xr{min}(R)}C_\xf{r}=\overline{K}$, and $T_\xf{r}$ denotes the integral closure of $R$ in $C_\xf{r}$. Again, the pullback $T_0=\{(\alpha_\xf{r})_{\xf{r}\in\xr{min}(R)}\in\prod_{\xf{r}\in\xr{min}(R)}T_\xf{r}\mid\varphi(\alpha_\xf{p}) =\psi(\alpha_\xf{q})\}$ is an aic of $R$ that is not isomorphic to $T$ itself, as it has only half the required idemps. We observe that the image of $R$ in $T$ is contained in $T_0$, as $\forall_{r\in R}\,\varphi(r)=r\text{ mod }\xf{m}=\psi(r)$. To see tightness of $T_0$ over $R$, let $0\ne t_0=(\alpha_\xf{r})_{\xf{r}}\in T_0$. Say $\alpha_\xf{r}\ne0$. This coordinate zeroes a monic $f\in R[X]$, for which we may assume $f(0)\notin\xf{r}$ (as in the second proof of the tightness of $T$ over $R$ in section \ref{sec:redux}). Take $c\in(\bigcap_{\,\xf{s}\in\xr{min}(R)-\{\xf{r}\}}\xf{s})-\xf{r}$. Then $cf(t_0)=0$ in $T_0$, hence $cf(0)\in (t_0T_0\cap R)-\xf{r}$, because (the image of) $R\subseteq T_0$. Ai closedness of $T_0$ follows as above, looking at just the $\xf{p}$-th and $\xf{q}$-th coordinates. But this contradicts (1).
\end{proof}
\end{thm}

\section{Using model theory}\label{sec:mt}

The concepts of tightness and algebraicity (rather than integrality) of ring extensions can be represented in first-order logic as the omission of two 1-types, as we will see in a moment. Based on this, we present a case for the conjecture that ``most'' rings in fact possess at least two non-isomorphic aics. For the benefit of the reader, a concise brush-up on the concepts is included a little further down.\\

We will adopt the terminology from Chang \& Keisler's \cite{CK}, except that, as in the formulation and the very name of the Omitting Types Theorem, by a ``1-type'' we shall mean just a consistent set $\Sigma$ of $\xc{L}$-formulas $\eta(x)$ that have at most one free variable: $x$, rather than a maximal consistent such set.\\

To start off, given a ring $R$, let $\xc{L}$ be the first-order language consisting of the function symbols $+$ and $\times$, plus for every $r\in R$ an individual constant $r^\bullet$. And let $\xc{T}$ be the first-order theory with axioms: the usual ones for commutative unital rings, plus the sentences $r_1^\bullet+r_2^\bullet=(r_1+r_2)^\bullet$ and $r_1^\bullet\times r_2^\bullet=(r_1r_2)^\bullet$ for all $r_1,r_2\in R$, plus $r_1^\bullet\ne r_2^\bullet$ for all $r_1,r_2\in R$ with $r_1\ne r_2$. The models of $\xc{T}$ are just the commutative rings $S$ that contain $R$ as a subring. We suppress the multiplication symbol $\times$ in formulas, as customary, and for the bullets $^\bullet$ we do the same.\\

For a sentence $\sigma$ of $\xc{L}$, i.e. a formula without free variables (variables not bound by $\forall$ or $\exists$), one writes $\xc{T}\vdash\sigma$ to say that $\sigma$ follows from (the axioms of) $\xc{T}$ under the derivation rules of first-order logic. If $\eta(x)$ is an $\xc{L}$-formula, $S$ a $\xc{T}$\!-model and $s\in S$, we say that $s$ \textit{realizes} $\eta(x)$, and write $S\models\eta(s)$, if $s$ satisfies the property expressed by $\eta(x)$ when interpreted in $S$. If such $S$ and $s$ exist, $\eta(x)$ is called \textit{consistent} with $\xc{T}$\!. This can be generalized to sets $\Sigma$ of formulas in a single variable $x$, called 1-\textit{types}. $s$ realizes $\Sigma$ in $S$ when $s$ realizes every formula $\eta(x)\in\Sigma$. And if such $S$ and $s$ exist, $\Sigma$ is said to be consistent with $\xc{T}$\!. When there is no such $s$, $S$ \textit{omits} $\Sigma$. Here and there, we enclose formulas in corners. This allows ``$=$'' to be used both as the formal equality symbol in formulas and for assignment statements in the metatheory.\\

As an example, let $\Sigma$ be the collection of the formulas $\varphi_{\vec{r}}(x)=\ulcorner\!r_nx^n+r_{n-1}x^{n-1}+\cdots+r_{0}=0\to r_n=0\urcorner$ for all $\vec{r}=(r_n,\cdots,r_{0})\in R^{n+1}$ of (arbitrary) length $n+1$, and take $R=\xb{Q}$. A $\xb{Q}$-algebra $S$ that omits $\Sigma$ is algebraic over $\xb{Q}$. For if $s$ is in $S$, there must be an $n$ and an $\vec{r}=(r_n,\cdots,r_0)\in\xb{Q}^{n+1}$ such that $\varphi_{\vec{r}}(s)$ fails in $S$. Then we have $r_ns^n+r_{n-1}s^{n-1}+\cdots+r_{0}=0$, but $r_n\ne0$. And $S$ realizes $\Sigma$ iff it contains an element that satisfies none of the equations, i.e., that is transcendental over $\xb{Q}$.\\

Next, we add to $\xc{T}$ enough sentences to make models $S$ ai closed. One sentence is needed for all (or infinitely many) polynomial degrees $n$. E.g., for the quadratic case it would be $\forall_x\forall_y\exists_u\exists_v(u+v=-x\wedge uv=y)$, so that every $Z^2+xZ+y$ will factor as $(Z-u)(Z-v)$ in $S[Z]$ for suitable $u,v\in S$.\\

We will denote the integral closure of $R$ in a $\xc{T}$\!-model $S$ by $S^\dagger$. By Lemma~\ref{lem:1}, the ring $S^\dagger$ is again ai closed.\\
 
For $S$ to be an aic of $R$, it needs to omit two 1-types of $\xc{L}$. The first one consists of one formula $\eta_{\vec{r}}(x)$ for every finite sequence $\vec{r}=(r_{n-1},\cdots,r_{0})\in R^n$ of any length $n$, saying that $x^n+r_{n-1}x^{n-1}+\cdots+r_{0}\ne0$. When $S$ does not realize this type, $S$ is integral over $R$. However, in cases like $R=\xb{Z}$, $\xc{T}$ \textit{locally realizes} this type, as the saying goes, since $\gamma(x)=\ulcorner2x=1\urcorner$ is a formula consistent with $\xc{T}$, and $\xc{T}\vdash\forall_x(\gamma(x)\to\eta_{\vec{r}}(x))$ for every $\vec{r}$, because $\frac{1}{2}$ (or any $s\in S$ with $2s=1$) cannot be a zero of any monic $f\in\xb{Z}[X]$, hence realizes the type in $S$. The formula $\gamma(x)$ describes what an element $x$ realizing the type may look like.\\

So instead, we will use the $\varphi_{\vec{r}}(x)$ defined above. Let $\Sigma_1$ be the set of these formulas, for all $n\in\xb{N}$ and $\vec{r}=(r_n,\cdots,r_{0})\in R^{n+1}$. Models $S$ that omit this type are algebraic over $R$. Instead of $S$ itself, the subring $S^\dagger$ will be ai closed and integral over $R$. $\xc{T}$ \textit{locally omits} $\Sigma_1$, that is, there is no $\xc{L}$-formula $\gamma(x)$ having at most $x$ free and consistent with $\xc{T}$ such that $\xc{T}\vdash\forall_x(\gamma(x)\to\varphi_{\vec{r}}(x))$ for all applicable finite sequences $\vec{r}$. For $\gamma(x)$ would then imply that $x$ is transcendental over $R$, and no single formula can of course achieve that on its own (except when $R=1$, where $\gamma(x)=\ulcorner0=1\urcorner$ achieves this tour de force).\\

For the definition of the second 1-type $\Sigma_2$, we first specialize to the case where $R=P=\prod_{a\in A}R_a$ is a finite product of non-trivial rings, as before, but we no longer require that the $R_a$ are domains, merely that they are connected. If $h=2^{|A|}$, then $R$ will have precisely $2^h$ idempotents (including the two trivial ones), say $d_1,\cdots,d_h$. The $e_a$ (for $a\in A$), introduced in Section \ref{sec:first}, are among them. All $\xc{T}$-models $S$ contain the fundamental system of orthogonal idempotents $\{e_a\mid a\in A\}$, so can be written as $\prod_{a\in A}S_a$, with $S_a=e_aS$. But these $S$ may of course have additional idempotents.\\

Let $\delta(x)$ abbreviate the $h$-fold conjunction $\bigwedge_{i=1}^hx\ne d_i$.\\

For $\vec{r}=(r_n,\cdots,r_{0})\in R^{n+1}$, put $\vartheta_{\vec{r}}(y)=\ulcorner\!y^n+r_{n-1}y^{n-1}+\cdots+r_{0}=0\urcorner$ (in which $r_n$ is not referenced), and for $n,m\in\xb{N}$, $\vec{r}\in R^{n+1}$ and $\vec{r}\,'\in R^{m+1}$, put $\psi_{\vec{r},\vec{r}\,'}(x)=\ulcorner\delta(x)\,\wedge\,\forall_y\,\forall_z\,((\vartheta_{\vec{r}}(y)\wedge\vartheta_{\vec{r}\,'}(z)\wedge xy=r_n\wedge(1-x)z=r_m')\to\bigwedge_{a\in A}(r_ne_a=0\vee r_m'e_a=0))\urcorner$, and let $\Sigma_2$ be the collection of all formulas $\psi_{\vec{r},\vec{r}\,'}(x)$ in $x$.\\

If $S$ omits $\Sigma_2$, every $s\in S$ not among the $d_i$ is (in particular) $\ne0$ and $1$, and there are $y,z\in S$ integral over $R$ plus an $a\in A$ such that $sy=r_n$ and $(1-s)z=r_m'$\! and $r_ne_a\ne0\ne r_m'e_a$ (for suitable $\vec{r}$ and $\vec{r}\,'$), meaning that if, in particular, $s\in S^\dagger$, both ideals $e_asS^\dagger$ and $e_a(1-s)S^\dagger$ have nonzero contraction to $R$. So the extension $S^\dagger/R$ is tight, and thus $S^\dagger$ is an aic of $R$.\\

If the $R_a$ are domains, then $\gamma(x)=\ulcorner\delta(x)\,\wedge\,x^2=x\!\urcorner$ implies $\psi_{\vec{r},\vec{r}\,'}(x)$ under $\xc{T}$ for all finite sequences $\vec{r}$ and $\vec{r}\,'$\! of this kind. For, arguing in $\xc{T}$\!, if $x$ is an idempotent that differs from all the $d_i$, and $y$ and $z$ (satisfy the monic equations involved and) have $xy=r_n$ and $(1-x)z=r_m'$\!, and there is an $a\in A$ with $r_ne_a\ne0$ and $r_m'e_a\ne0$, then $r_nr_m'\!=x(1-x)yz=0$. But since the $a$-th coordinates of $r_n$ and $r_m'$\! are non-zero, the fact that $r_nr_m'\ne0$ is known to $\xc{T}$, so that is a contradiction. This $\gamma(x)$ is consistent with $\xc{T}$\!, because $R\hookrightarrow R\times R$ via $r\mapsto(r,r)$ and $R\times R$ has extra idempotents: the $(e,e')$ with $e$ and $e'$\! different idempotents of $R$. So any ai closed $S\supseteq R\times R$ is a $\xc{T}$-model that realizes $\gamma(x)$. Hence $\xc{T}$ locally realizes $\Sigma_2$ when every $R_a$ is a domain.\\

And $\gamma(x)=\ulcorner\delta(x)\,\wedge\,\bigwedge_{(r,r')\in(R-\{0\})^2}\forall_y\,\forall_z\,(xy\ne r\,\vee\,(1-x)z\ne r')\urcorner$ does the same job when $R$ is a \textit{finite} ring $\ne1$. Indeed, $R\subsetneq R\times R$ and, for $s=(1,0)$, $(s)\cap R=0=(1-s)\cap R$, so $R\times R\models\gamma(s)$. So here $\xc{T}$ locally realizes $\Sigma_2$ too.\\

But for a general $R$, an $\xc{L}$-formula $\gamma(x)$ that locally realizes $\Sigma_2$ over $\xc{T}$ must be consistent with $\xc{T}$\!. I.e., there must be an extension $R\subseteq S$, with $S$ ai closed, and an $s\in S$, such that $S\models\gamma(s)$. And for every such $S$ the existence of such an $s$ needs to imply $S\models\psi_{\vec{r},\vec{r}\,'}(s)$ for \textit{all} pairs $(\vec{r},\vec{r}\,')$. If $r_n=0$ or $r_m'\!=0$, then $\bigwedge_{a\in A}(r_ne_a=0\vee r_m'e_a=0)$ is true, hence, by definition, it is implied by any formula $\gamma(x)$ whatsoever. And, if $\gamma(x)$ implies $\delta(x)\,\wedge\,x^2=x$, when $r_nr_m'\ne0$ we can argue as in the case of a finite product of domains.\\

However, since $\vartheta_{\vec{r}}(y)$ expresses that $y$ is integral with a specific equation over $R$, and the notion of being an element integral over $R$ with unknown equation cannot be formalized using a single formula, in case $r_n\ne r_nr_m'\!=0\ne r_m'$\! the only chance $\gamma(x)$ has, is to entail that $xy=r_n\wedge(1-x)z=r_m'$\! will be false for \textit{every} $y$ and $z$ in $S$, not just the integral ones. But that has proved to be a hard act - for this author at least, in view of the following remark.\\

\textit{Note.} Some prior versions of the paper claimed that also when $R$'s zero-divisors are easily surveyable, $\xc{T}$\! locally realizes $\Sigma_2$. But the proof that the formula $\gamma(x)$ used is consistent with $\xc{T}$\! was incorrect, and doesn't appear to be easily mended.\\

All in all, in light of the above considerations it seems quite plausible that for any infinite ring which is not a finite product of domains, the corresponding theory $\xc{T}$\! must locally omit $\Sigma_2$.\\

We now produce the (un)desired non-isomorphic aics, for rings $R$ for which $\xc{T}$\! \textit{does} locally omit $\Sigma_2$.

\begin{prop}\label{prop:8}Let $R$ be a finite product of connected rings, and put $\sigma_{\,\xr{extra}}=\ulcorner\,\exists_x(\delta(x)\,\wedge\,x^2=x)\urcorner$, a statement saying that extraneous idempotents exist (that are not already present in $R$ itself). If $\xc{T}_{\,\xr{extra}}=\xc{T}\cup\{\sigma_{\,\xr{extra}}\}$ and $\xc{T}_{\,\xr{nextra}}=\xc{T}\cup\{\neg\,\sigma_{\,\xr{extra}}\}$ locally omit $\Sigma_2$, then $R$ admits both aics with and aics without additional idempotents. Aics for $R$ are therefore \textit{not} unique.
\begin{proof}
We have already noticed that $\xc{T}$ locally omits $\Sigma_1$. So do the two extended theories. If $R$ is countable as a set, so is the language $\xc{L}$, hence, by the Extended Omitting Types Theorem (Th.~2.2.15 in \cite{CK}), under the given circumstances each of $\xc{T}_{\,\xr{extra}}$ and $\xc{T}_{\,\xr{nextra}}$ has a model omitting both types $\Sigma_1$ and $\Sigma_2$. The integral closures of $R$ in such models are aics of $R$, as we saw earlier. But idempotents are integral over $R$ and the $\xc{T}_{\,\xr{nextra}}$-model will have the same (finite number of) idempotents as $R$, while the $\xc{T}_{\,\xr{extra}}$-model has more. So the aics of $R$ that live in these two models can't be isomorphic.

For uncountable $R$, we can take $\kappa=|R|$, the cardinality of $R$, and instead use the $\kappa$-OTT (this is Th.~2.2.19 in Chang \& Keisler), by which every consistent theory $\xc{T}$ in a language $\xc{L}$ of power $\kappa$, which $\kappa$-omits an $n$-type $\Sigma$ of $\xc{L}$, has a model of cardinality $\leq\kappa$ that omits $\Sigma$. Here, $\xc{T}$ is said to $\kappa$\textit{-omit} $\Sigma$ if there is no set $\Gamma=\Gamma(x_1,\cdots,x_n)$ of formulas that have at most $x_1,\cdots,x_n$ free, with $|\Gamma|<\kappa$, and $\Gamma$ consistent with $\xc{T}$, and $\xc{T}\cup\Gamma(x_1,\cdots,x_n)\models\Sigma(x_1,\cdots,x_n)$. Now one cannot express transcendence over $R$ in fewer than $\kappa$ statements, and the same goes for \textit{looseness} (= untightness). The $\kappa$-OTT mentions only a single type $\Sigma$, while we use two. But this is easily overcome by letting $\Sigma$ be the $1$-type consisting of all disjunctions $\varphi_{\vec{r}\,''}(x)\vee\psi_{\vec{r},\vec{r}\,'}(x)$, where the $\psi_{\vec{r},\vec{r}\,'}$ are as above and $\vec{r}\,''$\! runs over all finite sequences in $R^{k+1}$, using a third (independent) natural number $k$. When a $\xc{T}$\!-model $S$ omits this $\Sigma$, for every $s$ in $S$ there are $n,m,k$ and vectors $\vec{r}$, $\vec{r}\,'$\!, $\vec{r}\,''$\! of lengths $n+1,m+1,k+1$, for which $\varphi_{\vec{r}\,''}(x)\vee\psi_{\vec{r},\vec{r}\,'}(x)$ fails for $x=s$, meaning that both disjuncts fail, meaning that etc. etc.
\end{proof}
\end{prop}

\end{document}